\documentstyle[12pt, amstex]{amsart}
\renewcommand{\pf}{\noindent{\em Proof: }}

\newcommand{\s}[1]{\mathcal{#1}}

\newcommand{\diag}{\operatorname{diag}}
\newcommand{\rank}{\operatorname{rank}}

\newcommand{\End}{\operatorname{End}}
\newcommand{\alg}{\operatorname{Alg}_{\mathbb C}}
\newcommand{\Span}{\operatorname{Span}}

\newtheorem{Thm}{Theorem}[section]
\newtheorem{Def}[Thm]{Definition} \newtheorem{Rem}[Thm]{Remark}
\newtheorem{Lem}[Thm]{Lemma} 
\newtheorem{Prop}[Thm]{Proposition}

\numberwithin{equation}{section}

\begin{document}

\title[Q-Heisenberg space]{Quantized Heisenberg Space}
\author[Jakobsen and Zhang]{Hans Plesner Jakobsen and Hechun
  Zhang$\vphantom{}^{1}$} 
\address{
  Mathematics Institute\\ Universitetsparken 5\\ DK--2100 Copenhagen {\O},
  Denmark} 
\address{Dept. of applied Math, Tsinghua
University\\ Beijing, 100084, P.R. China}\email{jakobsen@@math.ku.dk, hzhang@@math.tsinghua.edu.cn}

\date{\today}

\begin{abstract} We investigate the algebra $F_q(N)$ introduced by
  Faddeev, Reshetikhin and Takhadjian \cite{frt}. In case $q$ is a
  primitive root of unity the degree, the center, and the set of
  irreducible representations are found. The Poisson structure is
  determined and the De Concini-Kac-Procesi Conjecture is proved for
  this case. In the case of $q$ generic, the primitive ideals are
  described. A related algebra studied by Oh is also treated.
\end{abstract}
\maketitle

\footnotetext[1]{The second author
is partially supported by NSF of China}

\vskip 0.5cm
\medskip

\section{introduction}
\medskip
Among the  candidates for a quantized phase space is the
algebra $F_q(N)$ introduced by   Faddeev, Reshetikhin and Takhadjian
\cite{frt}. It is promising since it both contain generators
satisfying relations 
\begin{equation}
A_1B_1=qB_1A_1
\end{equation}
and generators satisfying 
\begin{equation}
A_2B_2-qB_2A_2=I.
\end{equation}
The higher dimensional analogues of the former are often known as ``the
quantum affine space'' (c.f. S. P. Smith, \cite{MR94g:17032}) whereas the
latter is among the candidates for ``the quantized Heisenberg algebra''.  We
shall give the precise defining relations of $F_q(N)$ below in
Section~\ref{2}. In \cite{frt} it was called the quantized function algebra of
Hermitian space - a name we hesitate to adopt. One reason is, as proved by
Korogodsky and Vaksman \cite{kv}, that the algebra is connected with
$SU_q(n,1)$ and hence not directly connected with a hermitian symmetric
space. Thus the name in the title.

The algebra $F_q(N)$ is very closely related to the high dimensional
quantum Heisenberg algebra (\cite{kv}) as well as to the so\--called
quantum symplec\-tic spa\-ce ${\s O}_q({\mathfrak sp}{\mathbb
C}^{2n})$. It is thus natural to expect that many features of those
algebras will show up in $F_q(N)$ as indeed they do. We shall return
to this connection when appropriate.

In the present article we analyze algebraic properties (center, degree) as
well as the representation theory of $F_q(N)$ when $q$ is a primitive root of
unity. The algebra fits into the general framework of Procesi and De Concini
(\cite{cp}). This first gives a tool for computing the center and the degree.
Moreover, and more importantly, it yields a Poisson structure which,
hypothetically is connected to the representation theory through the so-called
De Concini-Kac-Procesi Conjecture (\cite{kdp}). We prove that this conjecture
is true in the present setting.

Section~{\ref{2} gives the defining relations for the algebra and its
  associated quasipolynomial algebra. The relation to other algebras,
  especially the one considered by Oh (\cite{Oh}) is given.
  Section~\ref{3} describes the quite intricate arguments needed for
  obtaining explicitly the Poisson structure. The DKP Conjecture is
  then stated. Finally, the symplectic leaves are described, including
  their dimensions. In Section~\ref{4} the theory of De Concini and
  Procesi is briefly described and is then applied in Section~\ref{5}
  to finding the degree and more generally the canonical form of a
  certain number of quasipolynomial algebras. In Section~\ref{6}, the
  irreducible representations are found and the conjecture is
  verified. It is made precise how the algebra in \cite{Oh} also is
  covered by these results. Finally, in Section~\ref{7} we return
  briefly to the $q$ generic case. In \cite{Oh}, Oh described
  completely the set of primitive ideals for ${\s O}_q({\mathfrak
  sp}{\mathbb C}^{2n})$. We show briefly that in our case, the
  primitive ideals can be described in a similar way.

\medskip

\section{definitions and basic properties}

\label{2}
\medskip

The quantum function algebra of the quantum space ${\mathbb C}^N$ is an
associative algebra generated by
$z_0, z_1,\cdots,z_{N-1},z_0^*
z_1^*,\cdots,z_{N-1}^*$ subject to the following relations:
\begin{eqnarray}\label{jzrem} z_iz_j&=&q^{-1}z_jz_i\text{ for }i<j,\\\nonumber
      z_i^*z_j^*&=&qz_j^*z_i^*\text{ for }i<j,\\\nonumber
 z_iz_j^*&=&q^{-1}z_j^*z_i\text{ for }i\ne j,\textrm{ and }\\\nonumber
 z_iz_i^*-z_i^*z_i&=&(q^2-1)\sum_{k>i}z_kz_k^*,\end{eqnarray}
where $q\in{\mathbb C}^*$ is the quantum parameter.

\begin{Thm} The quantum function algebra $F_q(N)$ is an iterated Ore
extension and hence Noetherian and an integral domain with the same
Hilbert series of the commutative polynomial algebra in $2N$
variables.  There is a ${\mathbb C}$-basis for $F_q(N)$ consisting of
\begin{equation}{\s A}=\{z_0^{r_0}\cdots z_{N-1}^{r_{N-1}} {z^*_0}^{r^*_0}\cdots
{z^*_{N-1}}^{r_{N-1}^*}|r_i,r_i^*\in{\mathbb Z}_+\}.\end{equation}
\end{Thm}

\pf Observe that $F_q(N)$ is an iterated Ore extension
\begin{equation}F_q(N)={\mathbb C}[z_{N-1}][z_{N-1}^*,\sigma_{N-1}^*,\delta_{N-1}^*]
\cdots[z_0,\sigma_0,\delta_0] [z_0^*,\sigma_0^*,\delta_0^*]
\end{equation}
for automorphisms $\sigma_i,\sigma_i^* $ and $\sigma$-derivations
$\delta_i,\delta_i^*$ defined as follows:
\begin{eqnarray}\sigma_i:{\mathbb C}[z_{N-1},z_{N-1}^*,\cdots,z_{i+1},z_{i+1}^*]
&\longrightarrow & {\mathbb C}[z_{N-1},z_{N-1}^*,\cdots,z_{i+1},z_{i+1}^*],
\\\nonumber
\sigma_i(z_j)&=&qz_j\text{ for } j=i+1,\cdots,N-1,\\\nonumber
\sigma_i(z_j^*)&=&qz_j^*\text{ for }j=i+1,\cdots,N-1,\\\nonumber
\delta_i:{\mathbb C}[z_{N-1},z_{N-1}^*,\cdots,z_{i+1},z_{i+1}^*]
&\longrightarrow &{\mathbb C}[z_{N-1},z_{N-1}^*,\cdots,z_{i+1},z_{i+1}^*],
\\\nonumber
\delta_i(z_j)&=&\delta_i(z_j^*)=0\text{ for }j=i+1,\cdots, N-1,
\\\nonumber
\sigma_i^*:{\mathbb C}[z_{N-1},z_{N-1}^*,\cdots,z_{i+1},z_{i+1}^*,z_i]
&\longrightarrow & {\mathbb C}[z_{N-1},z_{N-1}^*,\cdots,z_{i+1},z_{i+1}^*,z_i],
\\\nonumber
\sigma_i^*(z_j)&=&q^{-1}z_j\text{ for } j=i+1,\cdots,N-1,\\\nonumber
\sigma_i^*(z_j^*)&=&q^{-1}z_j^*\text{ for }j=i+1,\cdots,N-1,\\\nonumber
\sigma_i^*(z_i)&=&z_i\\\nonumber
\delta_i^*:{\mathbb C}[z_{N-1},z_{N-1}^*,\cdots,z_{i+1},z_{i+1}^*,z_i]
&\longrightarrow &{\mathbb C}[z_{N-1},z_{N-1}^*,\cdots,z_{i+1},z_{i+1}^*,z_i],
\\\nonumber
\delta_i^*(z_j)=\delta_i^*(z_j^*)&=&0\text{ for }j=i+1,\cdots, N-1,\textrm{ and}
\\\nonumber
\delta_i^*(z_i)&=&(1-q^2)\sum_{k>i}z_kz_k^*.
\end{eqnarray}
This completes the proof. \qed

\medskip

Let $\overline{F_q(N)}$ be the associated quasipolynomial algebra of $F_q(N)$,
which is an associative algebra with generators $\bar{z}_0,\bar{z}_1,\cdots,
\bar{z}_{N-1}, \bar{z}_0^*,\bar{z}_1^*,\cdots, \bar{z}_{N-1}^*$ and the
following defining relations:

\begin{eqnarray} \bar{z}_i\bar{z}_j&=&q^{-1}\bar{z}_j\bar{z}_i\text{ for }i<j,\\\nonumber
      \bar{z}_i^*\bar{z}_j^*&=&q\bar{z}_j^*\bar{z}_i^*\text{ for }i<j,\\\nonumber
 \bar{z}_i\bar{z}_j^*&=&q^{-1}\bar{z}_j^*\bar{z}_i\text{ for }i\ne j,\textrm{ and}\\\nonumber
 \bar{z}_i\bar{z}_i^*&=&\bar{z}_i^*\bar{z}_i\text{ for all }i.\end{eqnarray}

\medskip

We finish this section by constructing an irreducible module for
$\overline{F_q(N)}$ of a ``high'' dimension:

\medskip

Let $q^m=1$. Let $\sigma, D\in \End({\mathbb  C}^m)$ be the operators which with respect to the standard basis
$v_0,v_1,\cdots,v_{m-1}$ are given as
\begin{equation}\sigma(v_{\underline{i}})=v_{{\underline{i+1}}},D(v_{\underline{i}})=q^iv_{\underline{i}}\quad \text{ for all }
 \quad {\underline{i}},\end{equation}where the index ${\underline{i}}$ in
$v_{\underline{i}}$ is $i\mod m$. We denote by $\sigma_i$ and similarly $D_i$ the operators $1\otimes
1\otimes\cdots\otimes\sigma\otimes 1\cdots\otimes 1$ and $1\otimes
1\otimes\cdots\otimes D\otimes 1\cdots\otimes 1$ on $({\mathbb  C}^m)^{N-1}$ with $\sigma$
and $D$ in the
{i}th positions for $i=1,\cdots,N-1$.
\begin{eqnarray}\bar{z}_i&=&D_{i}\sigma_{i+1}\sigma_{i+2}\cdots\sigma_{N-1}
\quad \text{ for }\quad i=1,\cdots,N-2,\\\nonumber
\bar{z}_0&=&D_{1}\sigma_{1}\sigma_{2}\cdots\sigma_{N-1},\\\nonumber
\bar{z}_{N-1}&=&D_{N-1},\\\nonumber
\bar{z}_{N-1}^*&=&D_{N-1},\\\nonumber
\bar{z}_{i}^*&=&D_{i}\sigma_{i+1}^{-1}\cdots\sigma_{N-1}^{-1}
  \text{ for } \quad i=1,2,\cdots,N-2,\textrm{ and }\\\nonumber
\bar{z}_{0}^*&=&D_1^{-1}\sigma_1^{-1}\cdots\sigma_{N-1}^{-1}.
\end{eqnarray}
Utilizing the relation $D\sigma=q\sigma D$, elementary computations now yield
that the above formulas define a representation of the algebra
$\overline{F_q(N)}$ in an $m^{N-1}$ dimensional space.  If the quantum
parameter $q$ is a primitive $m$th root of unity where $m$ is an odd positive
integer the representation is irreducible. In fact,
\begin{equation}\bar{z}_i\bar{z}_i^*=D_{i}^2\text{ for }i=1,2,\cdots,N-1,
\end{equation}
so we get every $D_i=(D_i^2)^{\frac{m+1}2}$ for $i=1,2,\cdots, N-1$. Then we
can get every $\sigma_i$ for all $i=1,2,\cdots,N-1$. Since $D,\sigma$ generate
$\End({\mathbb C}^m)$ the representation is irreducible.

\medskip

\begin{Rem}
The algebra considered in \cite{Oh} is given as follows (changing the notation in an inessential way): It is generated by $x_1, x_2,\cdots,x_{N},x_1^*,
x_2^*,\cdots,x_{N}^*$ subject to the following relations:
\begin{eqnarray} x_ix_j&=&q^{-1}x_jx_i\text{ for }i<j,\\\nonumber
      x_i^*x_j^*&=&qx_j^*x_i^*\text{ for }i<j,\\\nonumber
 x_ix_j^*&=&q^{-1}x_j^*x_i\text{ for }i\ne j,\textrm{ and }\\\nonumber x_i^*x_i-q^2
 x_ix_i^*&=&(q^2-1)\sum_{k<i}q^{k-i}x_kx_k^*.\end{eqnarray} Here, the
 factor $q^{k-i}$ may be gotten away with by absorbing $q^i$ into,
 say, $x_i$. After changing $q\mapsto q^{-1}$ and setting
 $z_i=x_{N-i}, z_i^*=x_{N-i}^*$ we get the following equivalent
 algebra generated by $z_0, z_1,\cdots,z_{N-1},z_0^*
 z_1^*,\cdots,z_{N-1}^*$ and with  relations:
\begin{eqnarray} z_iz_j&=&q^{-1}z_jz_i\text{ for }i<j,\\\nonumber
      z_i^*z_j^*&=&qz_j^*z_i^*\text{ for }i<j,\\\nonumber
 z_iz_j^*&=&qz_j^*z_i\text{ for }i\ne j,\textrm{ and }\\\nonumber
 z_iz_i^*-q^2 z_i^*z_i&=&(q^2-1)\sum_{k>i}z_kz_k^*.\end{eqnarray}
We return to this algebra in Remark~\ref{ohrem2}.\label{ohrem1}
\end{Rem}

\medskip

\section{The Poisson structure and its symplectic leaves}
\label{3}
\subsection{The DKP Conjecture}

Let $A$ be an associative algebra with generators $x_1,\cdots,x_n$ and defining relations
\begin{equation}\label{defrel}x_ix_j=q^{h_{ij}}x_jx_i+p_{ij},\text{ where }p_{ij}\in
{\mathbb C}[x_1,\cdots,x_{i-1}]\textrm{ if }i>j,\end{equation}
and where $(h_{ij})$ is an anti-symmetric integral matrix. Let $\varepsilon$ be a primitive $m$th root
of unity. Assume that every $x_i^m$ is central.
Let $\eta_1,\eta_2,\cdots,\eta_n$ be coordinate functions of ${\mathbb C}^n$
so that $\eta_i$ corresponds to $x_i^m$ for each $i=1,\cdots,n$. As explained
in \cite[p. 84-85]{cp}, there
is then a Poisson structure on ${\mathbb C}^n$ induced from the defining relations of the
algebra $A$ given by
\begin{equation}\{\eta_i,\eta_j\}=\left(\lim_{q\rightarrow \varepsilon}\frac{[x_i^m,x_j^m]}{m(q^m-1)}\right).
\end{equation}

\medskip

\begin{Rem}\label{w-r} Observe that up to a constant, the lowest order of $q^m-1$, when expanded around $\varepsilon$ is $q-\varepsilon$. Indeed, when $m$ is an odd prime and $\varepsilon\neq1$ is an $m$th root of unity, then $(q^m-1)=(q-1)(q-\varepsilon)\cdots (q-\varepsilon^{m-1})$. Moreover, 
\begin{equation}
\frac{q^m-1}{q-1}=1+q+\cdots+q^{m-1},
\end{equation}
hence
$(1-\varepsilon)(1-\varepsilon^2)\cdots(1-\varepsilon^{m-1})=m$. Hence,
$(q^m-1)\simeq\frac{m}{\varepsilon}(q-\varepsilon)$. It is also possible to let
$q\rightarrow 1$, but in this case, the Poisson structure is rather
degenerate (c.f. below).
\end{Rem}

A connected submanifold which is a maximal integral manifold of the
distribution defined by the Hamiltonian vector fields corresponding to
the Poisson structure, is called a {\bf symplectic leaf}. These define
a foliation of ${\mathbb C}^n$.

Let $\pi$ be an irreducible representation of the algebra $A$. Then there
exist a $p=(p_1,p_2,\cdots,p_n)$ in ${\mathbb C}^n$ such that $x_i^m=p_i$ on
$\pi$. Assume that ${\s O}_{\pi}$ is a symplectic leaf containing the point
$p$.

{\bf The DKP Conjecture:} For ``good'' $m$, \[\dim\pi=m^{\frac{1}{2}\dim{\s
    O}_{\pi}}.\]

\medskip

\begin{Rem}\label{dkprem}
An algebra with defining relations as in (\ref{defrel}) and where
$p_{ij}=0$ for all $i>j$ is called a quasipolynomial algebra. It is
easy to see that the DKP conjecture holds in this case provided that
$m$ is prime to all the numbers $h_{ij}$. Indeed, one may, without
loss of generality, assume that all the generators in an irreducible
module are invertible. But then the dimension of the irreducible
module is given as $m^{\frac12 r}$, where $r$ is the rank of the
matrix $\{h_{ij}\}$ (c.f. Section~\ref{4}) whereas the dimension of
the symplectic leaf is equal to the rank of the matrix
$\{b_{ij}=\{a_i,a_j\}\}$ where $b_{i,j}=h_{ij}a_ia_j$ and where for all $i$:
$a_i\equiv x_i^m$. Clearly the two matrices have the same rank.
\end{Rem}

\medskip

In the conjecture, the point only enters through the leaf to which it
belongs. To understand this, and also for further use, we need some
observations from \cite{cp}:

 Assume that we have a manifold $M$ and a
vector bundle $V$ of algebras with $1$ (i.e. $1$ and the
multiplication map are smooth sections). We identify the functions on
$M$ with the sections of $V$ which are multiples of $1$. Let $D$ be a
derivation of $V$, i.e. a derivation of the algebra of sections which
maps the algebra of functions on $M$ into itself. Let $X$ be the
corresponding vector field on $M$ and let $(t,p)\mapsto X_t(p)$ denote
the local $1$-parameter group determined by $X$.

\begin{Prop} For each point $p\in M$ there exists a neighborhood $U_p$ and a map
$(t,v_p)\mapsto\phi(t,v_p)$ defined for $|t|$ sufficiently small and
$v_p\in V_{\mid U_p}$ which  for each $t$ is a morphism of vector
bundles and which covers the $1$-parameter group generated by $X$. Indeed, for
each such $t$ and each $p\in U_p$, $\phi_t$ induces an algebra
isomorphism between $V_p$ and $V_{X_t(p)}$. \end{Prop}

Now, suppose $M$ is a Poisson manifold. Assume furthermore that the
Poisson structure lifts to $V$ i.e. that each local function $f$
induces a derivation on sections extending the given Hamiltonian
vector field $X_f$ determined by $f$. 

\begin{Prop}Under the above hypotheses, the fibers of $V$ over the points of
a given symplectic leaf $M$ are all isomorphic as algebras.\end{Prop}

In the present general set-up, $M={\mathbb C}^n$ and for each
$\underline{y}=(y_1,\cdots,y_n)\in {\mathbb C}^n$, the fiber
$V_{\underline{y}}$ of $V$ over $\underline{y}$ is given as 
\begin{equation}
V_{\underline{y}}=A/<x_i^m-y_i\ ;\; i=1,\cdots,n>.
\end{equation}
It was proved in \cite{cp} that the Poisson structure does lift in this situation.

\medskip

On a symplectic leaf $\s O$ the fibers of the algebra are isomorphic. Later
on, when studying the representations of the algebra $F_q(N)$ connected with
some leaf,  our strategy  will be to choose a good point $p\in{\s O}$.

\subsection{The Poisson structure defined by $F_q(N)$}

\medskip

Let $a_i,a_i*\leftrightarrow z_i^m,{z_i^*}^m$ for $I=0,1,\cdots,N-1$.

\begin{Prop}\label{post} The Poisson structure is given by
\begin{eqnarray*}
\{a_i,a_j\}&=&-a_ia_j \quad\text{ if $i< j$},\\
\{a_i^*,a_j^*\}&=&a_i^*a_j^* \quad\text{ if $i< j$},\\
\{a_i,a_j^*\}&=&-a_ia_j^* \quad\text{ if $i\neq  j$,\quad and}\\
\{a_i,a_i^*\}&=&2\sum_{k>i}a_ka_k^*. 
\end{eqnarray*}
\end{Prop}

In order to prove this proposition we now state and prove a number of Lemmas.

\begin{Lem}\label{s-lem}For any positive integer $s$
\begin{eqnarray}z_i(z_i^*)^s&=&(z_i^*)^sz_i+(q^{2s}-1)\sum_{k>i}z_kz_k^*,
(z_i^*)^{s-1}\textrm{ and}\\\nonumber
z_i^*z_i^s&=&z_i^sz_i^*+(q^{-2s}-1)\sum_{k>i}z_k^*z_kz_i^{s-1}.
\end{eqnarray}
In particular, if $q$ is an $m$th root of unity then the elements $z_i^m$,
$(z_i^*)^m$, and $z_{N-1}^a(z_{N-1}^*)^{m-a}$ are central for all
$i=0,1,\cdots, N-2$ and all $a=0,1,\cdots, m$.\end{Lem}

\pf This follows easily by induction. \qed

\medskip

\begin{Lem} Let \begin{eqnarray}\Omega_i&=&\sum_{k=i}^{N-1}z_kz_k^*,\\\nonumber\Omega_i^*&=&\sum_{k=i}^{N-1}z_k^*z_k, \text{ and}\\\nonumber
\Omega&=&\Omega_0.
\end{eqnarray}
We have the following formulas:
\begin{eqnarray}\Omega_iz_k&=&z_k\Omega_i\text{ if }
i\le k\le N-1,\\\nonumber
\Omega_iz_k&=&q^2z_k\Omega_i \text{ if }0\le k\le i-1,\\\nonumber
 \Omega_iz_k^*&=&z_k^*\Omega_i\text{ if }
i\le k\le N-1,\\\nonumber
\Omega_iz_k^*&=&q^{-2}z_k^*\Omega_i \text{ if }0\le k\le i-1,\\\nonumber
z_iz_i^*-z_i^*z_i&=&(q^2-1)\Omega_{i+1},\text{ and}
\\\nonumber
q^2z_iz_i^*-z_i^*z_i&=&(q^2-1)\Omega_i.\end{eqnarray}
In particular, if $q$ is an $m$th root of unity, then $\Omega$ is central.
\end{Lem}

\pf This follows directly from the defining relations. \qed

\medskip

Let $p_i=(q^{-2i}-1)$. For $i,j\in {\mathbb N}$ set $S_{i,j}=\{a\in {\mathbb N}\mid 1\leq a\leq j,\text{ and }a\neq i\}$. Let 
\begin{equation}
F_{i,j}(n)=\frac{p_i^{n-1}}{\prod_{a\in S_{i,j}}(p_i-p_a)}.
\end{equation}

\begin{Lem}\label{a-ff}
Let
\begin{equation}(z_iz_i^*)^n=\sum_{t=0}^n a_t(n)z_i^t(z_i^*)^t
\Omega_{i+1}^{n-t}.\end{equation}
Then 
\begin{equation}\label{a-f}
a_t(n)=\sum_{b=1}^t F_{b,t}(n).
\end{equation}
Moreover, since the right hand side makes sense for all $n$, we can define the
left hand side by this formula. Doing this, we get that $a_t(s)=0$ for
$s=1,\cdots,t-1$.
\end{Lem}

\pf We have
\begin{eqnarray}(z_iz_i^*)^{n+1}=\sum_{t=0}^n a_t(n)z_i^t(z_i^*)^t
\Omega_{i+1}^{n-t}z_iz_i^*\\\nonumber
=\sum_{t=0}^n a_t(n)z_i^t(z_i^*)^tz_iz_i^*
\Omega_{i+1}^{n-t}\\\nonumber
=\sum_{t=0}^na_t(n)z_i^t(z_i(z_i^*)t+(1-q^{2t})\Omega_{i+1}
(z_i^*)^{t-1})z_i^*\Omega_{i+1}^{n-t}\\\nonumber
=\sum_{t=0}^n a_t(n)z_i^{t+1}(z_i^*)^{t+1}
\Omega_{i+1}^{n-t} +\sum_{t=0}^na_t(n)(q^{-2t}-1)z_i^t(z_i^*)^t
\Omega_{i+1}^{n-t+1}.\end{eqnarray}
Hence, if we set $a_{-1}(n)=a_{n+1}(n)=0$, we get
\begin{equation}\label{rec1}\forall t=0,\cdots, n+1:\qquad
  a_t(n+1)=(q^{-2t}-1)a_t(n)+a_{t-1}(n). 
\end{equation}
In particular, 
\begin{equation}\label{rec2}\forall n: a_n(n)=1\quad\text{ and }\quad \forall n: a_0(n)=0.
\end{equation}

To  check that (\ref{rec1}) is satisfied by (\ref{a-f}), observe that
$p_t=(q^{-2t}-1)$ and that 
\begin{equation}
a_{t-1}(n)=\sum_{b=1}^{t-1}\left(\frac{p_b^{n-1}}{\prod_{a\in
      S_{b,t-1}}(p_b-p_a)}\right)\frac{p_b-p_t}{p_b-p_t}=\sum_{b=1}^{t-1}\left(\frac{p_b^{n-1}(p_b-p_t)}{\prod_{a\in
      S_{b,t}}(p_b-p_a)}\right). 
\end{equation}

 As for
(\ref{rec2}) as well as the remaining assertion, observe that 

\begin{equation}\sum_{b=1}^t \frac{p_b^{s-1}}{\prod_{a\in S_{b,t}}(p_b-p_a)}
\end{equation}
is the $(t,s)$th entry of the matrix product $A^{-1}\cdot A$ where $A$ is the
(Vandermonde) matrix

\begin{equation} 
A=\begin{pmatrix}1&p_1&p_1^2&\cdots&p_1^{t-1}\\1&p_2&p_2^2&\cdots&p_2^{t-1}\\\vdots&\vdots&\vdots&\cdots&\vdots\\1&p_t&p_t^2&\cdots&p_t^{t-1}\end{pmatrix}.
\end{equation}
\qed

\begin{Lem}\label{m-po}
\begin{equation}\Omega_i^m=\sum_{k=i}^{N-1}z_k^m{z_k^*}^m.
\end{equation}
\end{Lem}

\pf
Observe that
\begin{equation}\Omega_i^m=(z_iz_i^*+\Omega_{i+1})^m.
\end{equation}
 Since $z_iz_i^*$ and $\Omega_{i+1}$  commute we have
 \begin{equation}\label{b}\Omega_i^m=\sum_{n=0}^m\begin{pmatrix}
 m\\n\end{pmatrix} (z_iz_i^*)^n\Omega_{i+1}^{m-n}.
 \end{equation}
Therefore
\begin{equation}\Omega_i^m=\sum_{n=0}^m\begin{pmatrix}m\\n\end{pmatrix}
\sum_{t=0}^na_t(n)z_i^tz_i^{*t}\Omega_{i+1}^{m-t}.\end{equation}
The coefficient of $z_i^tz_i^{*t}\Omega_{i+1}^{m-t}$ is
\begin{equation} \sum_{n=t}^m\begin{pmatrix}m\\n\end{pmatrix}a_t(n).
\end{equation} 
By Lemma~\ref{a-ff} this may be written
\begin{equation} \sum_{n=1}^m\begin{pmatrix}m\\n\end{pmatrix}\sum_{b=1}^t \frac{p_b^{n-1}}{\prod_{a\in S_{i,j}}(p_b-p_a)}=\sum_{b=1}^t\sum_{n=1}^m\begin{pmatrix}m\\n\end{pmatrix} \frac{p_b^{n-1}}{\prod_{a\in S_{i,j}}(p_b-p_a)}.\end{equation}
The coefficient is clearly 1 if $t=m$. For all other values of $t$ it is zero
since then $p_b\neq0$ and thus
\begin{equation}\sum_{n=1}^m\begin{pmatrix}m\\n\end{pmatrix}p_b^{n-1}=\frac1{p_b}((p_b+1)^m-1)=0.
\end{equation} 
\qed

\medskip

If $f$ is a polynomial in $q,q^{-1}$, denote by $[f]$ the value
\begin{equation}
[f]=\lim_{q\rightarrow\varepsilon}\frac{f(q)}{m\cdot(q^m-1)}.
\end{equation}

\medskip

When computing $z_k^m{z_k^*}^m-{z_k^*}^mz_k^m$ we need more than the
above $\Omega^m$, namely we need to show that all other terms vanish
at the root of unity.

To ascertain this, consider more generally ( with $z\leftrightarrow z_k$ and $\Omega\leftrightarrow \Omega_{k+1}$.)

\begin{Lem}\label{l-s} Assume that $m$ is odd and that $q$ is a primitive $m$th root of unity. Let
\begin{equation}
z^i{z^*}^s={z^*}^sz^i+\sum_{j=1}^id_{i,j}(s)\Omega^j{z^*}^{s-j}z^{i-j}
\quad(i\leq s). \end{equation}
Then $[d_{m,j}(m)]=0$ for $j=1,\cdots,m-1$ and
$[d_{m,m}(m)]=2$.
\end{Lem}

\pf Let
\begin{equation}
d_{i,j}(s)=(q^{2s}-1)\cdots(q^{2(s+1-j)}-1)c_{i,j}(s).
\end{equation}

The recursion relation is: $c_{i+1,j}=q^{-2j}c_{i,j}+q^{-2(j-1)}c_{i,j-1}$ with
$c_{i,0}=1\;(=d_{i,0})$, and, for consistency,  $c_{i,j}=0$ if $j>i$. It follows that $c_{i,i}=q^{i-i^2}$. If we set
$f_{i,j}=q^{2(i-1)j}c_{i,j}$ the recursion relation becomes
\begin{equation}\label{recu}f_{i+1,j}=f_{i,j}+q^{2i}f_{i,j-1},\quad f_{i,0}=1,\text{ and }f_{i,i}=q^{i^2-i}.\end{equation}

The solution is

\begin{equation}f_{i,j}=\sum_{\delta=j-1}^{i-1}\cdots\sum_{\gamma=2}^{\delta-1}
\sum_{\alpha=1}^{\gamma-1}\sum_{\beta=0}^{\alpha-1}q^{2(\alpha+\beta+\gamma+\cdots+\delta)}.\end{equation}

Another way of expressing $f_{i,j}$ is as follows: Let ${\s Y}_{i,j}$ denote
the set of ``Young diagrams'' with $j$ rows with $n_k$ boxes in row $k$ and
such that $i-1\geq n_1>\cdots>n_j\geq0$. If $y\in {\s Y}_{i,j}$ let
$A(y)=n_1+\cdots+n_j$. Then 
\begin{equation}\label{area}
f_{i,j}=\sum_{y\in {\s Y}_{i,j}}q^{2A(y)}.
\end{equation} 

It is easy to verify that (\ref{area}) yields the solution to (\ref{recu}).

Now observe that the function $q^{n_1+n_2+\cdots+n_j}$ is invariant under
permutations. Also note that since $1+q+\cdots+q^{m-1}=0$, the sum of all $q^{n_1+n_2+\cdots+n_j}$ over the
set $P^+=\{(n_1,n_2,\cdots, n_j)\mid \forall k=1,\cdots,j: 0\leq n_k\leq m-1\}$
is  zero.

But, for the same reason, also the sum of $q^{n_1+n_2+\cdots+n_j}$ over sets of the form
$P_\omega^+=\newline \{(n_1,n_2,\cdots, n_j)\mid \forall k=1,\cdots,j: 0\leq n_k\leq
m-1\text{ and } n_{\omega_1}=\cdots n_{\omega_r}\}$ is zero.

It follows easily from this that $[d_{m,j}(m)]=0$ for
$j=1,\cdots,m-1$. Finally observe that since $m$ is odd (see also
Remark~\ref{w-r}),
\begin{equation}
[d_{m,m}(m)]=\lim_{q\mapsto\varepsilon}\frac{(q^{2m}-1)\cdots(q^2-1)\cdot q^{m^2-m}}{m\cdot(q^m-1)}=2.
\end{equation}
\qed 

\medskip

\noindent{\em Proof of Proposition~\ref{post}:} It should now be clear
from Lemma~\ref{a-ff}-\ref{l-s} (see also Remark~\ref{w-r}) that it only remains to consider $[z_i^m,z_j^m]$ where $z_iz_j=qz_jz_i$. Indeed, this is a prototype for the remaining cases. Since 
\begin{equation}
[z_i^m,z_j^m]=(q^{m^2}-1)z_j^mz_i^m,
\end{equation} 
this case follows since clearly $[q^{m^2}-1]=1$. \qed

\subsection{Symplectic leaves}

Consider a point 
\begin{equation}
\underline{a}=(a_{0},a_{1},\cdots ,a_{N-1},a_{N-1}^{*},\cdots ,a_{1}^{*},a_{0}^{*})\in 
{\mathbb C}^{2N}.
\end{equation}

Let $i_0$ denote the biggest $i$ such that $a_i\cdot a_i^*\neq0$. Suppose
that $a_i^*\neq0$ for some $i<i_0$. Let ${k_1}$ be
the biggest such number (below $i_0$). The
Hamiltonian vector field corresponding to $a_{k_1}$ has the form

\begin{equation}
\sum_j\left(-a_{k_1}a_j\frac{\partial}{\partial a_j}-a_{k_1}a_j^*\frac{\partial}{\partial a_j^*}\right)+2
a_{i_0}\cdot a_{i_0}^*\frac{\partial}{\partial a_{k_1}^*}
\end{equation}

The integral curves ($t$ is complex) of this have the form

\begin{eqnarray}
a_j(t)&=&a_j(0)e^{-a_{k_1}t}\\
a_j^*(t)&=&a_j^*(0)e^{-a_{k_1}t} \quad(j\neq {k_1})\\
a_{k_1}^*(t)&=&\left\{\begin{array}{l}a_{k_1}^*(0)-\frac{a_{i_0}(0)a_{i_0}^*(0)}{a_{k_1}}(e^{-2a_{k_1}t}-1)\quad
    \text{ for }a_{k_1}\neq
    0\\a_{k_1}^*(0)+2a_{i_0}(0)a_{i_0}^*(0)t \quad\text{ for }a_{k_1}= 0\end{array}\right.
\end{eqnarray}

So, we can flow to a point where $a_{k_1}^*=0$ unless we are in the exceptional case where

\begin{equation} a_{{k_1}}(0)a_{{k_1}}^*(0)+a_{i_0}(0)a_{i_0}^*(0)=0.
\end{equation}

\begin{figure}
-- d -- b -- a -- c -- $a^*$ -- $b^*$ -- 
\medskip

\begin{center}
\begin{tabular}{|c|c|c|c|c|c|c|}\hline\hline
$\{\cdot,\cdot\}$&$\{d,\cdot\}$&$\{b,\cdot\}$&$\{a,\cdot\}$&$\{c,\cdot\}$&$\{a^*,\cdot\}$&$\{b^*,\cdot\}$\\\hline\hline
$d$\vphantom{\large P}&$0$&$bd$&$ad$&$cd$&$a^*d$&$b^*d$\\\hline\hline
$b$&$-bd$&$0$&$ab$&$bc$&$ba^*$&$-2aa^*$\\\hline\hline
$a$&$-ad$&$-ab$&$0$&$ac$&$0$&$ab^*$\\\hline\hline
$c$&$-cd$&$-bc$&$-ac$&$0$&$a^*c$&$b^*c$\\\hline\hline
$a^*$&$-a^*d$&$-ba^*$&$0$&$-a^*c$&$0$&$a^*b^*$\\\hline\hline 
$b^*$&$-b^*d$&$2aa^*$&$-ab^*$&$-b^*c$&$-a^*b^*$&$0$\\\hline\hline
\end{tabular}
\end{center}
\smallskip
\caption{For convenience, we list the Poisson brackets for the
configuration above}
\end{figure}

Consider again the point 
\[
\underline{a}=(a_{0},a_{1},\cdots ,a_{N-1},a_{N-1}^{*},\cdots
,a_{1}^{*},a_{0}^{*})\in {\mathbb C}^{2N}.
\]
Let, as before, $i_{0}$ be the biggest index such that
$a_{i}a_{i}^{*}\neq 0$. Let $ r_{0}$ denote the number of non-zero
coordinates (starred or unstarred) having an index greater than that
$i_{0}$. If $a_{i}a_{i}^{*}=0$ for all $i$, let $r_0$ denote the
number of non -zero coordinates.  Let $i_{1}$ denote the biggest index
among those $i<i_{0}$ for which
$a_{i}a_{i}^{*}+a_{i_{0}}a_{i_{0}}^{*}=0$.  Let $i_{2}$ denote the
biggest index among those $i<i_{1}$ for which $
a_{i}a_{i}^{*}+a_{i_{1}}a_{i_{1}}^{*}+a_{i_{0}}a_{i_{0}}^{*}\neq 0$,
let $ r_{2}$ denote the number of non-zero coordinates $a_{i}$ and
$a_i^*$ with $ i_{2}<i<i_{1}$. Continuing like this we obtain a
sequence of indices $0\leq i_{s}<i_{s-1}<\cdots <i_{1}<i_{0}\leq
N-1$. And we set more generally $r_{2j}$ equal to the total number of
non-zero coordinates $a_{i}$ and $a_{i}^{*}$ with $i_{2j}<i<i_{2j-1}$.

Observe first that these quantities are invariant under Hamiltonian flow.
Also, we may, according to the preceding, use Hamiltonian flow to move to
a point in the same leaf where all coordinates $a_{i}\,$and $a_{i}^{*}$,
respectively, with $i_{2j+1}<i<i_{2j}$ are zero. We do that. Let us then
take a closer look at the $2N\times 2N$ matrix whose $ij$th entry is $
\left\{ x_{j},x_{i}\right\} $ where $x_{i}=a_{i-1}\,$for $i=1,\cdots ,N$,
and $x_{i}=a_{2N-i}^{*}\,$for $i=N+1,\cdots ,2N.$ The rank of this matrix is
(of course) equal to the dimension of the symplectic leaf. The coordinates $
a_{i}\,$and $a_{i}^{*}$ with $i_{2j+1}<i<i_{2j}$ now have non-zero Poisson
bracket with each other and zero with any other coordinate. Thus, they
contribute with a total of $2(i_{2j}-i_{2j+1}-1)$ to the rank.

Next observe that $\{a_{i_{1}},b\}$, where $b$ is some other coordinate
different from $a_{i_{0}}$ and $a_{i_{0}}^{*}$, is proportional to
$\{a_{i_{0}},b\}$.  The same observation holds for $a_{i_{1}}^{*}$ and
$a_{i_{0}}^{*}$ and even $a_{i_{0}}$ and $a_{i_{0}}^{*}$ have similar
Poisson brackets. Using this, it is easy to use row and column moves
to decouple the matrix into the direct sum of three submatrices - one
of rank 2, one involving the non-zero coordinates with $i>i_{0}$
together with $a_{i_{0}}^{*}$ and with a standard quasi-polynomial
coupling, and one involving the remaining coordinates, i.e. those
$a_i$ and $a_i^*$ with $i<i_1$.  The second summand involves $r_{0}+1$
points, and thus contributes with $2\cdot
\left[\frac{r_{0}+1}2\right]$ to the total rank. The last summand may
then be attacked with the same strategy, and after a finite number of
steps the result is obtained.

The total contribution to the rank from the first round is

\begin{equation}
2(i_0-i_1)+2\cdot \left[\frac{r_{0}+1}2\right].
\end{equation}

\medskip

If $s=2l$ let $i_{2l+1}=0$ and if $s=2l+1$ set $r_{2l+2}+1$ equal to the total
number of points having an index smaller than $i_{2l+1}$.

We thus get the following result:

\begin{Prop}\label{rank}Expressed in terms of the introduced data, the
  dimension of the symplectic leaf ${\s O}_{\underline{a}}$ containing
$\underline{a}$ is given by
\begin{equation*}
\dim{\s O}_{\underline{a}}=\left\{\begin{array}{l} 2\cdot \left[\frac{r_{0}}2\right]\text{ if
      }\forall i: a_{i}a_{i}^{*}=0.\\\sum_{k=0}^l \left(2(i_{2k}-i_{2k+1})+2\cdot
    \left[\frac{r_{2k}+1}2\right]\right) \textrm{ if } s=2l.\\\sum_{k=0}^{l} 2\left(i_{2k}-i_{2k+1}\right)+\sum_{k=0}^{l+1}\left(2\cdot
    \left[\frac{r_{2k}+1}2\right]\right) \textrm{ if } s=2l+1.\end{array}\right.
\end{equation*}
\end{Prop}

\section{The theory of De Concini and Procesi}
\label{4}

\medskip

The main tool used to compute the degrees and centers is the theory 
developed in \cite{cp} by De Concini and Procesi. Indeed, it is
straightforward to verify that $F_q(N)$ and $\overline{F_q(N)}$ satisfy the
hypotheses of that article and hence by \cite[6.4 Theorem]{cp} have the same
degree. 
 
Furthermore, there is a bijective correspondence between their centers by
means of the highest order term (c.f. \cite{jaz1}). 

\medskip 
 
Given an $n\times n$ skew-symmetric matrix $H=(h_{i,j})$ over ${\mathbb Z}$ one
constructs the twisted polynomial algebra ${\mathbb C}_H[x_1,x_2,\cdots,x_n]$
as follows: It is the algebra generated by elements $x_1,x_2,\cdots,x_n$ with
the following defining relations:
\begin{equation}x_ix_j=q^{h_{i,j}}x_jx_i\text{ for }i,j=1,2,\cdots,n.\end{equation} 
It can be viewed as an iterated twisted polynomial algebra with 
respect to any ordering of the indeterminates $x_i$. Given 
$a=(a_1,a_2,\cdots,a_n)\in{\mathbb Z}^n$ we write 
$x^a=x_1^{a_1}x_2^{a_2}\cdots x_n^{a_n}$. 
 
\medskip The degree and center of such an algebra is then given by \cite[7.1 Proposition]{cp}.

It is well known  that a skew-symmetric matrix over ${\mathbb Z}$ such as the 
matrix $H$ can be brought into a block diagonal form by an element 
$W\in SL({\mathbb Z})$. Specifically, there is a $W\in SL({\mathbb Z})$ and a sequence 
of $2\times2$ matrices 
$S(m_i)=\left(\begin{array}{cc}0&-m_i\\m_i&0\end{array}\right), i=1,\cdots,N,$ with 
$m_i\in {\mathbb Z}$ for each $i=1,\cdots,N$, such that   
\begin{equation} 
W\cdot H\cdot W^t=\left\{\begin{array}{l}{\text{\em Diag\;}}(S(m_1),\cdots, 
S(m_N),0)\text{ with $N=\frac{n^2-1}{2}$, if $n$ is odd}.\\{\text{\em Diag\;}}(S_1(m_1),\cdots, S(m_N))\text{  with $N=\frac{n^2}{2}$, if $n$ is even}.\end{array}\right. \label{diafo} 
\end{equation}

\begin{Def} We call the matrix $H$ the defining matrix. Any matrix of the form of the right-hand-side in (\ref{diafo}) will be called 
  a canonical form of $H$ and will occasionally be denoted by $J_H$. 
\end{Def}

\medskip

Thus, a canonical form of $H$ reduces the algebra to the tensor product of 
twisted Laurent polynomial algebras in two variables with commutation relation 
$xy=q^{r}yx$. By \cite[7.1 Proposition]{cp} it follows in particular that the 
degree of a twisted Laurent polynomial algebra in two variables is 
equal to $m/(m,r)$, where $(m,r)$ is the greatest common divisor of $m$ and 
$r$. 

\medskip

\section{ the center and the degree of the algebra $F_q(N)$}

\label{5}
\medskip

For use in determining the degree of $F_q(N)$ as well as for later, we now consider a number of quasipolynomial algebras.

\medskip

As usual, we let $\Omega_i=\sum_{k=i}^{N-1}z_kz_k^*$. Observe that
$\Omega_0,\Omega_1,\cdots,\Omega_{N-1}$ form a commutative family.

Consider $n$ points; $z_{i_1},\cdots,z_{i_a},\cdots,z_{i_n}$, with
$0\leq i_1<i_2<\cdots<i_n\leq N-1$ together with
$\Omega_{j_1},\Omega_{j_2},\cdots,\Omega_{j_r}$, with $0\leq
j_1<j_2<\cdots<j_r\leq N-1$. We wish to compute the degree of the
quasipolynomial algebra generated by these elements.

For $\ell=2,\cdots,r$, let $s_\ell$ denote the number of elements from
our family $z_{i_1},\cdots,$ $\cdots,z_{i_a}\cdots,z_{i_n}$ that have an index
$i_a$ satisfying $j_{\ell-1}\leq i_a<j_{\ell}$, let $s_1$ denote the
number of elements with an index $i_a$ satisfying $i_a<j_1$, and let
$s_{r+1}$ denote the number of elements with an index $i_a\geq
j_r$. To avoid redundancy, we will from now on {\em assume} that
$s_\ell\neq0$ for $\ell=1,2,\cdots,r$. (Indeed, if $s_\ell=0$ then
$\Omega_{j_\ell}$ can be removed from the algebra, c.f. below). If
$s_{r+1}\neq 0$ we indicate this with an $\downarrow$ and if $s_{r+1}=
0$ we indicate this with an $\uparrow$. More precisely, we denote the
algebras corresponding to these two cases by ${\s
L}_\downarrow(s_1,s_2,\cdots,s_{r},s_{r+1})$ and ${\s
L}_\uparrow(s_1,s_2,\cdots,s_{r})$, respectively. Observe that the
algebra is completely determined by the listed data. Finally, let ${\s
L}^T_{\downarrow\uparrow}={\s
L}_{\downarrow\uparrow}(\underbrace{1,1,\cdots,1}_T)$ for $T\in
{\mathbb N}$.

\begin{Prop}\label{n-prop}
  The non-trivial blocks in a canonical form of the defining matrix of
  ${\s L}^T_{\downarrow}$ consists of $\left[\frac{T-2}2\right]$
  blocks of $\begin{pmatrix}0&-4\\4&0\end{pmatrix}$,
  $\left[\frac{T}2\right]$ blocks of
  $\begin{pmatrix}0&-1\\1&0\end{pmatrix}$, and (only if $T$ is odd)
  one $\begin{pmatrix}0&-2\\2&0\end{pmatrix}$. In particular, the rank
  of ${\s L}^T_{\downarrow}$ is equal to $2T-2$. The non-trivial
  blocks in a canonical form of the defining matrix of ${\s
  L}^T_{\uparrow}$ consists of $\left[\frac{T}2\right]$ blocks of
  $\begin{pmatrix}0&-4\\4&0\end{pmatrix}$, $\left[\frac{T}2\right]$
  blocks of $\begin{pmatrix}0&-1\\1&0\end{pmatrix}$, and (only if $T$
  is odd) one $\begin{pmatrix}0&-2\\2&0\end{pmatrix}$. In particular,
  the rank of ${\s L}^T_{\downarrow}$ is equal to $2T$.
\end{Prop}

\pf When we in the sequel say ``by subtracting $a$ from $b$'' we mean
that we perform row and column operations corresponding to subtracting
the row and column of $a$ from the corresponding of $b$. The operation
of ``clearing'' is defined as follows: Suppose that we have a row
$r_1$ with a non-zero entry $w$ in the $c_1$th position and such that
the entries in the $c_1$th column are integer multiples of $w$.
Clearing then means subtracting appropriate multiples of this row from
the other rows, thereby creating a matrix whose $c_1$th column
consists of zeros except at the $r_1$th place. Of course, these
operations should be accompanied by corresponding column
operations. In our situations the row we start from will have just one
or two non-zero entries.

By subtracting $z_{T-1}$ from $z_T$ and clearing it is easy to see that one
can move from  ${\s L}^T_{\downarrow}$ to ${\s L}^{T-2}_{\uparrow}$ while
picking up one $\begin{pmatrix}0&-1\\1&0\end{pmatrix}$ in the process. This
reduces the problem to that of ${\s L}^T_{\uparrow}$. Using the same kind of
moves one can here move from ${\s L}^T_{\uparrow}$ to  ${\s
  L}^{T-2}_{\uparrow}$ while picking up both a
$\begin{pmatrix}0&-1\\1&0\end{pmatrix}$ and a
$\begin{pmatrix}0&-4\\4&0\end{pmatrix}$. The result then follows by induction
by paying special attention to the configuration of ${\s
  L}^1_{\uparrow}$. Indeed, the latter is just
$\begin{pmatrix}0&-2\\2&0\end{pmatrix}$. \qed

\medskip

For later use, observe that

 \begin{Lem}\label{5.5}
Let $M_x=\sum_{x\geq i>j\geq 1}(E_{i,j}-E_{j,i})$. Let $S_1=\begin{pmatrix}0&-1\\1&0\end{pmatrix}$. A canonical form for
$M_x$ is $\diag(\underbrace{S_1,\cdots,S_1}_{\frac x2})$ if $x$ is even, and
$\diag(\underbrace{S_1,\cdots,S_1}_{\frac {x-1}2},0)$ if $x$ is odd.  
\end{Lem}

\pf Starting by subtracting the second row from the first, this follows easily
by simple row and column operations.  \qed.

\medskip

A canonical form of the defining matrix of the associated
quasipolynomial algebra follows from Proposition~\ref{n-prop} via the
following, easily
established result. 

\begin{Lem}
A canonical form of the defining matrix of the associated
quasipolynomial algebra is equal to that of ${\s L}^N_{\downarrow}$. If $m$ is
odd, the degree of $F_q(N)$ is equal to $m^{N-1}$. 
\end{Lem}

\medskip

\begin{Prop}
  The rank of the defining matrix of ${\s
    L}_\downarrow(s_1,s_2,\cdots,s_{r},s_{r+1})$ is
  $R_\downarrow=2\left[\frac{s_{r+1}}2\right]+2\left[\frac{s_r+1}2\right]+\cdots+ 2\left[\frac{s_1+1}2\right]$.
  Let $x$ denote the number of {\em odd} $s_i$. The defining matrix of
  \newline ${\s
  L}_\downarrow(s_1,s_2,\cdots,s_{r},s_{r+1})$ is equivalent to 
\begin{equation}
\sum_{i=1}^{r+1}\left[\frac{s_{i}}2\right]\text{ copies of
      }\begin{pmatrix}0&-1\\1&0\end{pmatrix}\oplus \left\{\begin{array}{l}{\s
          L}^x_{\uparrow}\;(s_{r+1}\text{ even}\;).\\{\s L}^x_\downarrow\;(s_{r+1}\text{ odd}\;).\end{array}\right.
\end{equation}

The rank of the
defining matrix of ${\s L}_\uparrow(s_1,s_2,\cdots,s_{r})$ is
$R_\uparrow=2\left[\frac{s_r+1}2\right]+\cdots+ 2\left[\frac{s_1+1}2\right]$.  Let $y$ denote the
number of {\em odd} $s_i$ ($i=1,\cdots,r$). The defining matrix of ${\s
  L}_\uparrow(s_1,s_2,\cdots,s_{r})$ is equivalent to 
\begin{equation}
\sum_{i=1}^{r}\left[\frac{s_{i}}2\right]\text{ copies of
      }\begin{pmatrix}0&-1\\1&0\end{pmatrix}\oplus {\s L}^y_{\uparrow}.
\end{equation}
(This result remains true even if $s_\ell=0$ for one or more $\ell=1,2,\cdots,r$.)\label{deg}
\end{Prop}

\pf First of all, if one is just interested in the rank (or the cases
$m$ odd), one may replace the $\Omega$'s by
$\Omega_{j_u}\Omega_{j_{u+1}}^{-1}$ for each pair $u=1,\cdots,r-1$
together with $\Omega_{j_1}$ and then subtract the resulting rows and
columns from everything to decompose the matrix into a direct sum of
matrices of the form of Lemma~\ref{5.5}.

In the general case, let us assume that we have an algebra ${\s
L}_{\downarrow\uparrow}(\cdots, s_k,\cdots)$ with $s_k\geq2$. Let
$j_{k-1}\leq i-1<i<j_k$. By subtracting $z_{i-1}$ from $z_i$ and
clearing, one sees that one may move from ${\s
L}_{\downarrow\uparrow}(\cdots, s_k,\cdots)$ to ${\s
L}_{\downarrow\uparrow}(\cdots, s_k-2,\cdots)$ while picking up a
$\begin{pmatrix}0&-1\\1&0\end{pmatrix}$. The result follows easily
from this. \qed

\medskip 

We shall be content to treat the center in the case where $m$ is odd:

\begin{Prop}
Let $q$ be a primitive $m$th root of 1 with $m$ odd. The center of $F_q(N)$ is
generated by the following elements
\begin{eqnarray}
&\Omega, z_i^m \text{ and } {z_i^*}^m \text{ with }i=0.1.\cdots,N-1,
\text{ and }\\\nonumber &(z_{N-1})^a(z_{N-1}^*)^b\text{ with }a+b=m\text{ and }
a,b\in \{1,2,\cdots, m-1\}.
\end{eqnarray}
\end{Prop}

\pf The analogous result for the quasipolynomial algebra follows immediately
from the construction of a canonical form. But by comparing the highest
order terms, it is easy to see that the centers of $F_q(N)$ and
$\overline{F_q(N)}$, respectively, have the same magnitude. \qed

\medskip

\section{the irreducible representations of the algebra $F_q(N)$}

\label{6}

\medskip

In this section, the De Concini-Kac-Procesi Conjecture is proved to be true
for the algebra $F_q(N)$.

If for some element $w$ in some algebra, $uw=q^{\alpha_u}wu$ for all
generators $u$, then $w$ is said to be covariant. Observe than if in an
irreducible module, $w$ is covariant, then either $w=0$ or $w$ is
invertible.

Consider an irreducible module $V$. Let $z_i^m=a_i$ and ${z_i^*}^m=a_i^m$ for
all $i=0,1,\cdots, N-1$.  Construct from a point
\[\underline{a}=
(a_{0},a_{1},\cdots ,a_{N-1},a_{N-1}^{*},\cdots ,a_{1}^{*},a_{0}^{*})\in
{{\mathbb C}}^{2N}
\]
the sequence of indices $0\leq i_{s}<i_{s-1}<\cdots <i_{1}<i_{0}\leq N-1$ as
previously. Assume furthermore that we, possibly after using Hamiltonian flow,
are in a situation where  all coordinates $a_{i}\,$and $a_{i}^{*}$,
respectively, with $i_{2j+1}<i<i_{2j}$ are zero.

Let us first observe (c.f. (\ref{m-po})) that 
\[\forall i:\Omega^m_i=\sum_{j>i}a_ja_j^*.\]

Moreover, each $\Omega_i$ is covariant, hence if $\Omega_i^m=0$ in $V$, then
also $\Omega_i=0$ on $V$.

Let us then begin by considering points $z_i, z_i^*$ with $i>i_0$. It follows
that at most one of $a_i,a_i^*$ is non-zero and that $z_i$ commutes with
$z_i^*$. Suppose $a_i\neq0$. Then $z_i$ is invertible and $z_i^*$ is nilpotent
and ($q$-)commutes with everything. Hence $z_i^*=0$. Of course, if $a_i=0$
then also $z_i=0$.

The same kind of reasoning gives that $z_{i_0}$ and $z_{i_0}^*$ are
two commuting invertible operators. Here we send
$z_{i_0},\Omega_{i_0}$ to a set $\s A$ and discard $z^*_{i_0}$.

The operators $z_i$ with $i_1<i<i_0$ form a $q$-commuting family of nilpotent
operators as do the operators $z_i^*$ with the same restriction on the index.
Moreover, there is a common null space for these operators $z_i$ which is
invariant under all other operators except the corresponding $z_i^*$. Finally
observe that for any $i_1<i<i_0$, 
\[z_iz_i^*-z_i^*z_i=(q^2-i)z_{i_0}z_{i_0}^*\]
and the right hand side is invertible. Hence the operators $z_i,z_i^*$ in this
region are non-zero.

Let us now move to the operators $z_{i_1},z_{i_1}^*$. By assumption and the
above remarks concerning $\Omega_i$ it follows that
\[z_{i_1}z_{i_1}^* + z_{i_0}z_{i_0}^*=0.\]
The defining relations easily give that
$z_{i_1}z_{i_1}^*=q^{-2}z_{i_1}^*z_{i_1}$. Since all operators clearly are
invertible we may solve for $z_{i_1}^*$. The treatment of this operator is thus
completed and we will not consider it any more.

Moving in the direction of decreasing indices on the way to $i_2$ we
may pick up some terms $z_i$ or $z_j^*$ but never both
$z_k,z_k^*$. These are sent to $\s A$ and so are
$z_{i_2},\Omega_{i_2}$, but $z^*_{i_2}$ is discarded. After this, the
picture  repeats itself periodically. As a result, it
is clear that we obtain the following:

\begin{itemize} 
\item A family $N^+$ of non-zero $q$-commuting nilpotent operators $z_i$ with $i_{2j+1}<i<i_{2j}$.
\item A family $N^-$ of non-zero $q$-commuting nilpotent operators $z_i^*$ with
  $i_{2j+1}<i<i_{2j}$.
\item A family ${\s A}\subseteq \alg\{z_0,\cdots,z_{N-1},z_{N-1}^*,\cdots,z_0^*\}$ of
  ($q$-) commuting invertible operators.
\end{itemize}

\smallskip

These operators satisfy

\smallskip

\begin{enumerate} 
\item The operators from $N^+$, ${\s A}$, and $N^-$ generate all operators in the
  representation. 
\item For each $i_{2j+1}<i<i_{2j}:\;z_iz_i^*-z_i^*z_i=z_{i_{2j}}z_{i_{2j}}^*$.
\item ${\s A}$ leaves the common null spaces of $N^+$ and  $N^-$, respectively,
  invariant.
\item The only pairs of commuting operators from ${\s A}$ are
  ($z_{i_{2j}},z_{i_{2j}}^*$). 
\end{enumerate}

\medskip

Let $V_0$ denote the common null space of $N^+$ and let $\pi$ denote the
representation of ${\s A}$ on $V_0$. Then, clearly,

\begin{itemize}
\item $\pi$ is irreducible.
\item $V=N^-V_0$.
\end{itemize}

\medskip

By \cite[\S 7]{cp}, the representation $\pi$ is essentially unique.
Furthermore, if we set $N^-=\{w_1^*,\cdots,w_s^*\}$ then any element $v\in V$
can be written as a sum of elements
\begin{equation}{w_1^*}^{d_1}\cdots{w_s^*}^{d_s}v_d\label{ba}\end{equation}
where $\forall j=1,\cdots,s: 0\leq d_j\leq m-1$ and $v_d\in V_0$. Using the
action of $N^+$ it is easy to see that if  $v_d$ is taken from a fixed basis
of $V_0$ then the set of all elements (\ref{ba}) constitutes a basis of $V$. Observe that
$s=(i_0-i_1-1)+\cdots+(i_{2j}-i_{2j-1}-1)+\cdots$.

Conversely, given an irreducible ${\s A}$ representation $\pi$ on a
space $V_0$ promote this to a a representation of the algebra
$\alg\{{\s A},N^+\}$ by letting $N^+$ act trivially. Then it is
easy to see that the induced representation
\[\alg\{N^-,{\s A},N^+\}\otimes_{\alg\{{\s A},N^+\}} V_0\]
is irreducible. 

\medskip

\begin{Thm} Let ${\underline{a}}$ and $V$ be as above. Let $m$ be an odd integer. Then the DKP Conjecture
    is true in this case. Specifically,
\begin{equation}
\dim(V)=m^{\frac12 \dim {\s O}_{\underline{a}}}.
\end{equation} 
\end{Thm}

\pf A good deal of the proof has been done in the paragraphs proceeding the
statement of the theorem. What remains is to determine $\dim(V_0)$. This is
done by the results in Section~\ref{6} and by translating the notation of
Section~\ref{4} (and implicitly, the present) to that situation. 

Let us first consider the case where $\forall i:a_ia_i^*=0$. Our algebra is
then a quasipolynomial algebra and (possibly after a renumbering of e.g. the
$z_i^*$) the dimension is given by Lemma~\ref{5.5}.

In the remaining cases it suffices, due to the assumptions on $m$ to
use the rank of the defining matrix as given by
Proposition~\ref{deg}. Indeed, due to the way the algebra is
constructed, $s_{r+1}=r_0+1$ and thus we are always in the
``$\downarrow$ case''. If we are in the $s=2\ell$ case of
Section~\ref{4} then we get the following translation between the
remaining notation: $r=\ell+1$, $\forall
i=1,\cdots,r-1:s_{r+1-i}=2+r_{2i}$, and $s_1=0$. In this case the
contribution from $N^-$ to the dimension is
\begin{equation}
(i_0-i_1-1)+\cdots+(i_{2\ell-2}-i_{2\ell-1}-1)+(i_{2\ell}-i_{2\ell+1}),
\end{equation}
where, moreover, $i_{2\ell+1}=0$.  Thus, the two contributions to $\dim(V)$ do
not match up precisely with the two contributions to $\dim{\s
  O}_{\underline{a}}$. But the  discrepancies clearly cancel. 

In the remaining case, $s=2\ell+1$, the only differences are that
$s_1=r_{2\ell+2}+2$ and that the contribution to the dimension from $N^-$ now
is
\begin{equation}
(i_0-i_1-1)+\cdots+(i_{2\ell-2}-i_{2\ell-1}-1)+(i_{2\ell}-i_{2\ell+1}-1).
\end{equation}
The claim then follows as above. \qed

\medskip

\begin{Rem}\label{ohrem2}
As mentioned in the introduction, the algebra in \cite{Oh} is very
closely related to ours. Specifically, when this algebra is given in a
form as in Remark~\ref{ohrem1} it is easy to see how the
representation theory and symplectic leaves of this algebra are
connected with the ones given in the previous sections. Specifically,
consider an irreducible representation of the algebra in
Remark~\ref{ohrem1}. Without loss of generality we may assume that
$z_{N-1}^m\neq0$ or ${z_{N-1}^*}^m\neq0$. Since both
$z_{N-1},z_{N-1}^*$ are covariant we may in fact assume that at least
one of them is non-zero (hence invertible). The two cases are
equivalent, so let us just assume that $z_{N-1}\neq0$. Let $w_i=z_i$
and $w_i^*=z_i^*z_{N-1}^{2}$ for all $i=0,\cdots,N-1$. The resulting
relations are then,  in terms of $w_0,
w_1,\cdots,w_{N-1},w_0^* w_1^*,\cdots,w_{N-1}^*$, the following:
\begin{eqnarray}\label{modif} w_iw_j&=&q^{-1}w_jw_i\textrm{ for }i<j,\\\nonumber
      w_i^*w_j^*&=&qw_j^*w_i^*\text{ for }i<j<N-1,\\\nonumber
 w_iw_j^*&=&q^{-1}w_j^*w_i\text{ for }i\ne j\textrm{ and }i\neq
 N-1,\\\nonumber w_iw_i^*- w_i^*w_i&=&(q^2-1)\sum_{k>i}w_kw_k^*\textrm{
 for } i=0,\cdots,N-2,\\\nonumber
 w_{N-1}w^*_j&=&qw_j^*w_{N-1},\\\nonumber
 w_i^*w^*_{N-1}&=&q^3w^*_{N-1}w^*_i \textrm{ for }i\neq N-1,\textrm{ and}\\\nonumber
 w_{N-1}w_{N-1}^*&=&q^2 w_{N-1}^*w_{N-1}.\\\end{eqnarray}
In particular, if we let $\Omega_i=\sum_{k\geq i}w_kw_k^*$ (as before), then
\begin{eqnarray}
\forall i,j: \Omega_i\Omega_j&=&\Omega_j\Omega_i,\\\nonumber
\forall i: w_{N-1}\Omega_i&=&q^2\Omega_iw_{N-1},\textrm{ and }\\\nonumber
\forall i: w_{N-1}^*\Omega_i&=&q^{-2}\Omega_iw^*_{N-1}.
\end{eqnarray}
This maybe looks rather complicated, but it is in fact not. In dealing
with the representation theory, one can proceed with the same strategy
as for the previous case. Even more so, after a small adjustment, the
same results can be used. Since the Poisson structure (connected with
(\ref{modif})) behaves in the same way as far as the terms
$\{a_i,a^*_i\}$ are concerned - the fact that we have
$\{a_{N-1},a^*_{N-1}\}=2a_{N-1}a^*_{N-1}$ is of no importance here -
again allows one to move to a good point where the representation
theory breaks up in two parts. The part corresponding to $N^{\pm}$ is
identical whereas for the algebra corresponding to $\s A$ we get again
a family consisting of some elements $w_i,\Omega_j$ with $i,j\leq N-2$
together with $w_{N-1},w_{N-1}w_{N-1}^*$. Here we may set
$\Omega_{N-1}=w_{N-1}w_{N-1}^*$ whereas to deal with the remaining
generator we add a new variable $w_{-1}=w_{N-1}^{-1}$. (Formally, we
should then change the numbering by one unit on all variables, but we
shall not do so.) We can then take over directly the results of
Section~\ref{5}. All algebras will be of the form ${\s
L}_{\uparrow}$. The addition of the extra variable of course means a
shift in various quantities. For instance, since $\rank {\s
L}_{\uparrow}(\underbrace{1,1,\cdots,1}_N)=2N$, the degree is $m^N$.

\smallskip

As for the dimensions of the symplectic leaves, observe (this we could
also have done for the original algebra) that at a good point, the
Poisson structure decouples corresponding to the algebra breaking into
$N^{\pm}$ and ${\s A}$. The part corresponding to $N^{\pm}$ is
completely trivial to deal with. For the remaining part we have a
collection of functions $a_i,a_j^*$. Among the members of this family
are $a_{N-1},a_{N-1}^*$ and $a_{N-1}a_{N-1}^*\neq 0$ at the
point. After this there may be variables $a_{i_1}, a_{i_1}^*$ such
that $a_{i_1}a_{i_1}^*+a_{N-1}a_{N-1}^*= 0$, and so on, together
possibly with some extra functions $a_k$ and/or $a_r^*$. The crucial
point is that the functions $a_{i_s}$, $a_k$, and $a_r^*$ correspond
exactly to the variables $w_{i_s}$, $w_k$, and $w_r^*$ in $\s A$. At
the same time, for $j=i_0,i_2,\cdots$, the product $a_{j} a_{j}^*$ corresponds
exactly to $\Omega_j$. Moreover, the original functions can be
recovered from these. But this turns the problem into comparing
dimensions of representations and symplectic leaves for
quasipolynomial algebras, and here the result is well known (see
Remark~\ref{dkprem}).

\smallskip

In particular, {\bf The DKP conjecture holds for this algebra}.
\end{Rem}

\medskip

\section{Annihilators}

\label{7}

In this section we study annihilators of simple modules in the case where $q$
is generic. If $u$ is invertible,
$\textrm{Ad}(u):w\mapsto uwu^{-1}$ denotes as usual the adjoint
  representation of $u$. If $w$ is covariant we shall occasionally call
  $q^{\alpha_u}$ (or simply ${\alpha_u}$) the weight of $w$ w.r.t. $u$. 

We proceed by observing

\begin{Lem}\label{inv}
  Let ${\s C}$ be a quasipolynomial algebra with center $Z_{\s C}$ and $V_1$
  an irreducible module in which all generators are invertible. Let $\xi: Z_{\s
    C}\mapsto {\mathbb C}$ be the corresponding central character. Let
  $Z^\xi_{\s C}$ denote the subalgebra of $Z_{\s C}$ generated by
  $\{c-\xi(c)\mid c\in Z_{\s C} \}$. Then
\begin{equation}
{\textrm{Ann}}_{\s C}=F_q(N)\cdot Z^\xi_{\s C}.
\end{equation} 
\end{Lem}

\pf Let ${\s C}=\{w_1,\cdots,w_r\}$ and assume $w_iw_j=q^{a_{i,j}}w_jw_i$ all
$i,j=1,\cdots,r$. (It is a standing assumption that the coefficients $a_{i,j}$
are all integers.) The moves which brings the skew-symmetric matrix
$\{a_{i,j}\}$ into block diagonal form may simply be interpreted as a series
of replacements of the generators by ${\em monomials}$ in the generators and
their inverses; $w_i\mapsto u_i=w_1^{\alpha_{1,i}}\cdots w_r^{\alpha_{r,i}}$
for all $i=1,\cdots,r$ and with each $\alpha_{i,j}\in {\mathbb Z}$. Here we may
assume that $u_1u_2=q^{b_{1,2}}u_2u_1, \cdots,
u_{2s+1}u_{2s+2}=q^{b_{2s+1,2s+2}}u_{2s+2}u_{2s+1}$ ($2s+2\leq r$) and with
all other relations trivial (commutative). Now suppose that $u$ is an element
in the annihilator. Clearly the new generators are invertible and each
${\textrm{Ad}}(u_i)$ leaves the annihilator invariant. Hence we may assume
that $u$ has a fixed weight w.r.t. each generator $u_i$. But then evidently
$u$ must be a polynomial in the generators $u_{2s+3}, \cdots, u_r$. These
generators clearly generate the center and thus the proof is complete. \qed

\medskip

From now on, let $V$ be an irreducible module of the  algebra $F_q(N)$ and
let $\textrm{Ann}$ denote the annihilator of $V$ in $F_q(N)$.

Observe that if for some index $i$, $\Omega_{i+1}\neq 0$, then $z_i,z_i^*$ do
not commute and are not identically 0. If $\Omega_{i}=0$ we get that
$[z_i,z_i^*]= (\Omega_i-z_iz_i^*)$  and hence $z_iz_i^*=q^{-2}z_i^*z_i$.

\begin{Lem}\label{lem1} Suppose that $\Omega_{i+1}\neq 0$, $\Omega_{i}\neq 0$, and that  
\[\sum_\alpha z_i^\alpha p_\alpha +\sum_\beta {z_i^*}^\beta
r_\beta\in\textrm{Ann}\] where the $p_\alpha,r_\beta$ may contain powers of
$\Omega_i$ but otherwise only contain generators with index $j\neq i$. Then
each $z_i^\alpha p_\alpha$ and ${z_i^*}^\beta r_\beta$ belongs to ${\textrm
  Ann}$. Moreover, we may assume that each $p_\alpha$ and each $r_\beta$ is
covariant w.r.t. $\Omega_1$ and $\Omega_{i-1}$. 
\end{Lem}

\pf Notice that the element considered is the most general element containing
$z_i,z_i^*$. Observe that ${\textrm{Ad}}(\Omega_i)$ and
${\textrm{Ad}}(\Omega_{i+1})$ are identical on $p_\alpha,r_\beta$. Hence the
claims follows by weight considerations.   \qed

\medskip

We now introduce several subalgebras. Let ${\s B}$ denote the set of
generators $z_i,z_i^*$ for which $\Omega_{i+1}\neq 0$, $\Omega_{i}\neq 0$, and
let ${\s B}_0$ denote the subset of ${\s B}$ consisting of generators
$z_i,z_i^*$ such that one of them is not injective. Let ${\mathfrak b}$ and
${\mathfrak b}_0$ denote the corresponding sets of indices $i$ of the elements
in ${\s B}$ and ${\s B}_0$, respectively.  Let ${\s D}$ denote the
subalgebra generated by the elements of an index $i\notin {\mathfrak b}_0$
together with the operators $\Omega_i$ with $i\in{\mathfrak b}_0$. 
Let ${\s C}_1$ denote the subalgebra generated by the elements of an index
$i\notin {\mathfrak b}$ together with the operators $\Omega_i$ with
$i\in{\mathfrak b}$. The latter is a quasipolynomial algebra and hence its
generators are either invertible or identically zero. Let finally $\s C$
denote the algebra generated by the invertible elements in
$z_0,z_1,\cdots,z_1^*,z_0^*$ having an index $i\notin {\mathfrak b}$
together with the operators $\Omega_i$ with $i\in{\mathfrak b}$. Thus, ${\s
  C}$ is a  subalgebra of ${\s C}_1$.

We now prove some lemmas about the operators in these algebras.

\smallskip

\begin{Lem} If $i\in {\mathfrak b}_0$ then at most one of
  the pair $z_i,z_i^*$ is not invertible. 
\end{Lem}

\pf Let us for simplicity assume that $\forall i\in {\mathfrak b}_0: z_i$ not
injective. Let $V_0$ denote the common null space. This is invariant under all
generators whose index is not in ${\mathfrak b}_0$. Hence the space $V$ is of
the form
\[V=\Span \{{z_{i_1}^*}^{j_1}\cdots {z_{i_d}^*}^{j_d}v_{j_1,\cdots,j_d}\mid
i_1,\cdots,i_d\in {\mathfrak b}_0, {\textrm{ and }}v_{j_1,\cdots,j_d}\in V_0\}.\] 

Suppose that, say, ${z^*_{i_1}}$ annihilates a sum $\omega$ of such elements
and that at least one of the summands is non-zero. Let $k$ be the biggest
power of ${z^*_{i_1}}$ occurring in the expression (and assume that $k>0$).
Then $z_{i_1}^k\omega$ also is a sum in which at least one summand is
non-trivial.  Moreover, since $z_{i_1}^k{z_{i_1}^*}^k=\Omega_{i+1}^k +
(\cdot)z_{i_1}$ this leads to a new sum with fewer and still non-trivial summands
(recall that $\Omega_{i+1}$ is invertible and covariant). \qed

\medskip

 Actually, the former proof also yields 

\begin{Lem} $V_0$ is an irreducible ${\s D}$ module.
\end{Lem}

\medskip

\begin{Lem}\label{lem2}  ${\textrm{Ann}}$ is generated by those elements in the annihilator
  ${\textrm{Ann}}_{\s D}$ of $V_0$ in ${\s D}$ that are covariant with respect
  to the operators $z_i,z_i^*$ with $i\in{\mathfrak b}_0$.
\end{Lem}

\pf Let $i\in {\mathfrak b}_0$ and assume as before that $z_i$ is not
injective. Assume that ${z_i^*}^\beta r_\beta\in {\textrm{Ann}}$. Then by
injectivity of $z_i^*$, $r_\beta\in{\textrm{Ann}}$. In particular, $z_i
r_\beta\in {\textrm{Ann}}$. This case then becomes a special case of the
following: Suppose that $z_{i}^\ell p\in{\textrm{Ann}}$.  Then
$z_i^kp{z_i^*}^k\in {\textrm{Ann}}$ for all $k=\ell+1,\cdots$. Notice that if
$p$ is a sum of monomials then each monomial is covariant w.r.t. $z_i$. Define
$(p)_k$ by $z_i^kp=(p)_kz_i^k$. Observe that $p{z_i^*}^k={z_i^*}^k(p)_k$. Then
$(p)_k(z_i^k{z_i^*}^k-{z_i^*}^kz_i^k)\in{\textrm{Ann}}\subseteq {\textrm{Ann}}(V_0)$.  Since
$z_i^k{z_i^*}^k-{z_i^*}^kz_i^k=f(q)\Omega_{i+1}^k$ on $V_0$ with $f(q)\neq0$
(c.f. Lemma~\ref{s-lem}), each $(p)_k=0$. Hence we may assume that $p$
annihilates $V_0$ and has a fixed weight w.r.t. $z_i$ (it already has a fixed
weight w.r.t. $\Omega_i$ and $\Omega_{i+1}$). Hence, since
$z_iz_i^*=\Omega_i-\Omega_{i+1}$, it has a fixed weight w.r.t. $z_i^*$.
According to Lemma~\ref{lem1} we are now done with this index. Moving through
${\mathfrak b}_0$ the conclusion is easily reached.  \qed

\medskip

\begin{Lem}\label{lem3}  ${\textrm{Ann}}_{\s D}$ is generated by those elements in
  ${\textrm{Ann}}_{\s D}\cap {\s C}_1$ that are covariant with respect
  to the generators of ${\s D}$.
\end{Lem}

\pf This follows again from Lemma~\ref{lem1} as in the proof of Lemma~\ref{lem2}. \qed

\medskip

\begin{Prop}\label{la-pro}
The annihilator ${\textrm{Ann}}$ is generated by those generators $z_i$ and
$z_i^*$ of $F_q(N)$ that have an index $i\notin{\mathfrak b}$ and for which
$z_i$ and/or $z_i^*$ is identically zero on $V$, together with $\Omega-c$ for
some $c\in {\mathbb C}$. 
\end{Prop}

\pf By Lemma~\ref{lem2} and Lemma~\ref{lem3} an element $a$ in ${\textrm{Ann}}$
may be assumed to belong to ${\s C}_1$ and to have a specific weight with
respect to all generators (c.f. the proof of Lemma~\ref{inv}). Thus, if $a$ is a sum of monomials in the
generators, e.g. $a=a_1+a_2$, it may be assumed that $a_1$ and $a_2$ have the
same weight. By covariance,  $a_1$ and $a_2$ are either identically zero
or invertible. We may for now assume that they are both invertible. But then
$a_2\cdot a_1^{-1}$ is in the center of $F_q(N)$ and hence $a_2=a_1\cdot w$
for some $w$ in the center of $F_q(N)$. \qed

\medskip

We finish by some observations concerning the converse of
Proposition~\ref{la-pro}.

Let us then consider an irreducible representation $\hat V_0$ of ${\s
  C}_1$ (or, equivalently, an irreducible representation of ${\s C}$ in which
all generators are invertible). Let $i\in {\mathfrak b}\setminus{\mathfrak b}_0$. We want to have that 
\[\Span\{z_i^rv_r,{z_i^*}^sv_s\mid v_r,v_s\in \hat V_0\}\]
is an irreducible subspace for (ultimately) ${\s D}$, that $z_i,z_i^*$ are
injective, etc.

Utilizing the fact that $\Omega_i$ commutes with $z_i,z_i^*$ whereas
their weight w.r.t. $\Omega_{i+1}$ is $q^{\mp2}$, it follows that if
$z_i({z_i^*}^sv_s)= 0$ or $z_i^*({z_i}^rv_r)\neq 0$ for some $r,s$
then $\Omega_{i+1}= q^{n} \Omega_i$ for some $n\in\mathbb Z$. We will
from now on assume that $\forall n\in{\mathbb Z}: \Omega_{i+1}\neq
q^{n} \Omega_i$.

Elements of the form 
\[z_i^rv_r+z_i^{r-1}v_{r-1}+\cdots+v_0\]
are by definition non-zero (polynomials) and hence non-zero if just one
$v_i\neq0$. By applying $z_i$ (and $z_i^*$) appropriately, it follows that any
invariant subspace must contain such an element. Let $r$ be minimal
($r\geq1$). Observe that $\Omega_{i+1}=c^*\Omega_i$ for some non-zero complex
constant $c^*$.  Applying $\Omega_i$ and $\Omega_{i+1}$ to the element we must
still be in the invariant subspace. But then there is an element of the form
$z_i^\ell v_\ell$ in the subspace, and then an element $v_0$ in the
subspace. Hence, $\hat V_0$ is contained in any invariant subspace.

Observe that if $\Omega_{r+1}= q^{n} \Omega_i$ for some $n\in \mathbb Z$,
then $z_i$ and $z_i^*$ satisfy a covariance relation (which generically is not
homogeneous w.r.t. $q$).

But it follows easily that under the assumptions, $\Omega_i,\Omega_{i+1}$ are invertible  on the space  
\[z_i^rv_r+z_i^{r-1}v_{r-1}+\cdots+v_0+\cdots+{z_i^*}^sv_s.\]

\smallskip

Indeed, we obtain the following

\begin{Prop}
  Suppose that $\pi$ is an irreducible representation of ${\s C}$ for which
  $\forall i\in {\mathfrak b},\forall n\in {\mathbb Z}: \Omega_i\neq
  q^n\Omega_{i+1}$. Then there is an irreducible representation of $F_q(N)$
  for which the annihilator is generated by ${\s C}_1\setminus{\s C}$ together
  with $\Omega-c$. The complex constant $c$ is determined by $\pi$ if $0\in
  {\mathfrak b}$. Otherwise, $c=0$.
\end{Prop}

\medskip

Notice that this description is very similar to Oh's description by means of
{\em admissible sets}.

\end{document}